\documentclass[11pt, dvipsnames]{amsart}
\usepackage{amsthm}
\usepackage{amssymb}
\usepackage{amsmath}
\usepackage[margin=3.9cm]{geometry}

\usepackage[english]{babel}
\usepackage{tikz}
\usepackage[shortlabels]{enumitem}
\usepackage{url}
\usepackage{hyperref}
\usepackage{microtype} 
\usepackage[hyphenbreaks]{breakurl} 

\theoremstyle{plain}
\newtheorem{thm}{Theorem}[section]
\newtheorem{prop}[thm]{Proposition}
\newtheorem{lemma}[thm]{Lemma}
\newtheorem{cor}[thm]{Corollary}
\newtheorem*{mainthm}{Main Theorem}

\theoremstyle{definition}
\newtheorem{df}[thm]{Definition}
\newtheorem{remark}[thm]{Remark}
\newtheorem{example}[thm]{Example}
\newtheorem*{alg}{Algorithm}
\makeatletter
\newcommand{\customlabel}[2]{%
   \protected@write \@auxout {}{\string \newlabel {#1}{{#2}{\thepage}{#2}{#1}{}} }%
   \hypertarget{#1}{#2}%
}
\makeatother

\def\Q{\mathbb{Q}}
\def\Z{\mathbb{Z}}
\def\F{\mathbb{F}}

\DeclareMathOperator{\ICM}{ICM}
\DeclareMathOperator{\Pic}{Pic}
\DeclareMathOperator{\Hom}{Hom}
\DeclareMathOperator{\Tr}{Tr}

\newcommand{\frB}{{\mathfrak B}}
\newcommand{\frf}{{\mathfrak f}}

\newcommand{\p}{{\mathfrak p}}
\newcommand{\q}{{\mathfrak q}}
\newcommand{\frP}{{\mathfrak P}}
\newcommand{\cI}{{\mathcal I}}
\newcommand{\cJ}{{\mathcal J}}
\newcommand{\cK}{{\mathcal K}}

\newcommand{\cS}{{\mathcal S}}
\newcommand{\cO}{{\mathcal O}}

\newcommand{\cV}{{\mathcal V}}
\newcommand{\cW}{{\mathcal W}}

\renewcommand{\bar}{\overline}
\newcommand{\vphi}{{\varphi}}
\newcommand{\set}[1]{\left\lbrace#1\right\rbrace }
\newcommand{\typeover}[2]{\mathrm{type}_{#1}({#2})} 

\usepackage{color}

\title[Local isomorphism classes of fractional ideals]{Local isomorphism classes of fractional ideals of orders in \'etale algebras}
\date{\today}

\author{Stefano Marseglia}
\address{Mathematical Institute, Utrecht University, P.O. Box 80010, 3508 TA, Utrecht, The Netherlands}
\email{s.marseglia@uu.nl}



\begin{document}

\begin{abstract}
    We study the local isomorphism classes, also known as genera or weak equivalence classes, of fractional ideals of orders in \'etale algebras.
    We provide a classification in terms of linear algebra objects over residue fields.
    As a by-product, we obtain a recursive algorithm to compute representatives of the classes, which vastly outperforms previously known methods.
\end{abstract}

\maketitle


\section{Introduction}\
\label{sec:intro}
Let~$R$ be an order in a number field~$K$. 
Two fractional~$R$-ideals~$I$ and~$J$ are called \emph{isomorphic} if they are so as~$R$-modules, and they are said to be \emph{locally isomorphic} if for every maximal ideal~$\p$ of~$R$ we have~$I_\p\simeq J_\p$ as~$R_\p$-modules, where~$R_\p$ is the completion of~$R$ at~$\p$, $I_\p=I\otimes_R R_\p$ and $J_\p=J\otimes_R R_\p$.
Recall that a fractional~$R$-ideal is invertible if and only if it is locally principal, that is, locally isomorphic to~$R$. 
If~$R$ is the maximal order~$\cO$ of~$K$ then the set of isomorphism classes of fractional $R$-ideals is a finite commutative group, the class group of~$K$.
One of the characterizations of~$\cO$ is that every fractional~$\cO$-ideal is invertible.
As soon as we replace~$\cO$ with a non-maximal order~$R$ in~$K$, there will be more than one local isomorphism class of fractional~$R$-ideals.

The notion of local isomorphism class has assumed different names in the previous literature.
In \cite{DadeTausskyZas}, the authors used the expression \emph{weak equivalence}.
We will use the same terminology.
Important is also the notion of \emph{genus} of a fractional ideal.
More precisely, two fractional ideals~$I$ and~$J$ are in the same genus if for every rational prime~$p$ we have~$I\otimes_\Z \Z_p \simeq J\otimes_\Z \Z_p~$ as~$R\otimes_\Z \Z_p~$-modules.
On the one hand, as pointed out in \cite[Sec.~5]{LevyWiegand85}, two fractional ideals are weakly equivalent if and only if they are in the same genus.
On the other hand, the notion of genus can easily be generalized to non-commutative orders.
Accounts of known results can be found, for example, in \cite{Reiner03}, \cite{Roggenkamp70I}, \cite{Roggenkamp70II}, \cite{Guralnick84} and \cite{Guralnick87}.

In this paper we study the set~$W(R)$ of weak equivalence classes of fractional~$R$-ideals, which inherits the structure of a commutative monoid from the ideal multiplication.
As explained in detail in Section~\ref{sec:running_time_comparison}, the results contained in \cite{DadeTausskyZas} can be turned into an algorithm to compute~$W(R)$.
New results about the structure of $W(R)$ have been given by the author in \cite{MarsegliaICM} and \cite{2022arXiv2MarCMType}, resulting in faster algorithms.
In this paper, building on the Main Theorem stated below, we give a new algorithm to compute~$W(R)$, which, in general, vastly outperforms the previous ones.

The new method builds on two preliminary observations.
The first is that we have a decomposition
\[ W(R) = \bigsqcup_S \bar W(S), \]
where~$S$ runs over the overorders of~$R$ and~$\bar W(S)$ denotes the subset of~$W(R)$ consisting of the classes represented by fractional ideals~$I$ with~$(I:I)=S$.
The second observation is that we can reduce the computation of~$\bar W(R)$ to the computation of the weak equivalence classes of a finite number of overorders of~$R$, each having at most one non-invertible maximal ideal. 
See Theorem~\ref{thm:split}.
Combining the two observations, we are reduced to computing~$\bar W(R)$ for an order~$R$ with a unique non-invertible maximal ideal.
We are ready to state our main result, which is an immediate consequence of Theorem~\ref{thm:wk_classes_only_linear_algebra} and Proposition~\ref{prop:U_one_sing} presented later in the text.
As above, denote by~$R_\p$ the completion of~$R$ at a maximal ideal~$\p$.
Let~$k$ be residue field of~$\p$ and~$T=(\p:\p)$ the multiplicator ring of~$\p$.
Set~$T_\p = T\otimes_R R_\p$ and consider the $k$-algebra~$A=T/\p$. 
\begin{mainthm}
    Assume that~$\p$ is the unique non-invertible maximal ideal of~$R$.
    Denote by~$J_1,\ldots,J_n$ the representatives of~$\bar{W}(T)$.
    For each~$i=1,\ldots,n$, let~$\cV_i$ denote the set of sub-$k$-vector spaces~$\cI$ of~$J_i/\p J_i$
    satisfying
    \[ \set{a \in A : a \cI \subseteq \cI } = k\quad\text{and}\quad\cI A = J_i/\p J_i.\]
    Then the group~$U=T^\times_\p/R^\times_\p$ acts freely on each~$\cV_i$, and we have a bijection
    \[ \bar{W}(R) \longleftrightarrow \bigsqcup_{i=1}^n \frac{\cV_i}{U}. \]
\end{mainthm}
The Main Theorem, together with the results of Section~\ref{sec:split}, gives rise to a recursive algorithm~\ref{alg:WRbar} to compute~$\bar W(R)$ for an order~$R$ with any number of non-invertible maximal ideals.
This, in turn, gives us an algorithm~\ref{alg:WR} to compute the monoid~$W(R)$ of weak equivalence classes of~$R$.

These algorithms can be used to compute the monoid of isomorphism classes of fractional~$R$-ideals, which we denote by~$\ICM(R)$.
This is achieved by combining~$\bar W(S)$ and the group of invertible ideal classes~$\Pic(S)$, for every overorder~$S$ of~$R$, as explained in Section~\ref{sec:example_app}.

In the same section, we discuss a second application of the Main Theorem, namely computing isomorphism classes of abelian varieties over finite fields.
Indeed, the isomorphism classes of abelian varieties in isogeny classes satisfying certain assumptions are in bijection with the ideal class monoid of an order uniquely determined by the isogeny class.
Exploiting this connection, we compute the isomorphism classes of hundreds of thousands of isogeny classes of abelian varieties over finite fields.
The algorithms presented in this paper have been crucial to complete this computation, since previous versions were not able to finish for certain isogeny classes.
In a joint effort with Edgar Costa, Taylor Dupuy, David Roe, Christelle Vincent and Mackenzie West, we uploaded (a subset of) the result of this computation on the LMFDB \cite{lmfdb}\footnote{currently available through \url{https://abvar.lmfdb.xyz/Variety/Abelian/Fq/}}, together with other interesting arithmetic information about the varieties, like their endomorphism rings and polarizations.
See Section~\ref{sec:example_app} for more details.

The paper is structured as follows. 
In Section~\ref{sec:prelim}, we introduce the notation and recall some basic properties of orders and their fractional ideals.
We will not restrict ourselves to the number field case, but we will consider orders in \'etale algebras over the fraction field of an arbitrary Dedekind domain.  
In Section~\ref{sec:equiv}, we recall the definition of weak equivalence and introduce other notions of local equivalence.
In Section~\ref{sec:split}, we explain how to reduce the computation of the monoid of weak equivalence classes for a generic order to the case of an order with at most one non-invertible maximal ideal.
In Section~\ref{sec:mult_ring}, we discuss how to control the multiplicator ring of the extension of fractional ideals to overorders.
Here, we use that orders come equipped with a trace form.
These results are used in Section~\ref{sec:wk_ext} to show how these extensions with controlled multiplicator ring can be used to find representatives of each weak equivalence class.
Section~\ref{sec:wk_lin_alg} contains the main contribution of the paper, where we combine the results of Sections~\ref{sec:split} and~\ref{sec:wk_ext} to explain how to compute the weak equivalence classes of an order using only linear algebra over the residue fields of the non-invertible primes.
These results are turned into the algorithms~\ref{alg:WR} and~\ref{alg:WRbar} which are presented and discussed in Section~\ref{sec:alg}.
The running time of the algorithms is discussed in Section~\ref{sec:running_time_comparison}, together with a comparison with previously known algorithms.
Finally, Section~\ref{sec:example_app} contains examples and an explanation of how to apply the algorithms to compute ideal classes and abelian varieties over finite fields.

\subsection*{Acknowledgments}
The author is supported by NWO grant VI.Veni.202.107.
The author thanks Jonas Bergstr\"om, Valentijn Karemaker and John Voight for comments on a preliminary version of the paper.

\section{Preliminaries and notation}\
\label{sec:prelim}
Let~$Z$ be a Dedekind domain\footnote{for us, fields are not Dedekind domains. In particular, Dedekind domains have Krull dimension~$1$ and are infinite} and~$Q$ its fraction field.
Let~$K$ be an \emph{\'etale~$Q$-algebra}, that is, a finite product of separable finite field extensions of~$Q$.
A \emph{$Z$-lattice} in~$K$ is a finitely generated sub-$Z$-module of~$K$ which contains a~$Q$-basis of~$K$.
Consider the trace form~$\Tr=\Tr_{K/Q}$ defined by associating to an element~$a$ in~$K$ the trace of the matrix representing multiplication by~$a$ in~$K$ with respect to any~$Q$-basis of~$K$.
For every~$Z$-lattice~$I$, define the \emph{trace dual} as 
\[ I^t = \set{ a\in K :\Tr(aI)\subseteq Z }. \]
We have~$(I^t)^t = I$.

Given two~$Z$-lattices~$I$ and~$J$ in~$K$, define their \emph{colon ideal} as
\[(I:J) = \set{a\in K : aJ \subseteq I},\]
and, if~$J\subseteq I$, the \emph{index}~$[I:J]$ as the product of the invariant factors of the torsion~$Z$-module~$I/J$.

A \emph{$Z$-order}~$R$ in~$K$ is a subring which is also a~$Z$-lattice.
If the underlying Dedekind domain is clear from the context we will simply say \emph{order}.
All orders in~$K$ are contained in a unique order~$\cO$ which is called the \emph{maximal order} of~$K$.
If~$R$ and~$T$ are orders such that~$R\subseteq T$ then we say that~$T$ is an \emph{overorder} of~$R$.

A \emph{fractional~$R$-ideal} is a sub-$R$-module of~$K$ which is also a~$Z$-lattice in~$K$.
Given any~$Z$-lattice~$I$ in~$K$, then~$(I:I)$ is an order in~$K$, called the \emph{multiplicator ring} of~$I$.
Moreover, if~$I$ and~$J$ are fractional~$R$-ideals, then~$I^t$ and~$(I:J)$ are as well,
and we have natural isomorphisms
\[ \Hom_Z(I,Z) \simeq I^t \text{ and } \Hom_R(J,I) \simeq (I:J). \]

Colon and trace dual are related as explained in the following well known lemma; see for example~\cite[Lem.~15.6.12]{Voight21}.
\begin{lemma}\label{lemma:colon_dual}
    Let~$I$ and~$J$ be~$Z$-lattices in~$K$. Then 
    \[ I^t  J = (I:J)^t. \]
\end{lemma}

Let $\p$ be a maximal ideal of~$R$.
We denote by~$R_\p$ the completion of~$R$ at~$\p$ and by~$K_\p$ its total ring of quotients.

For a prime~$p$ of~$Z$, let~$Z_p$ be the completion of~$Z$ at~$p$ and~$Q_p$ the fraction field of~$Z_p$.
Note that~$K_\p$ is an \'etale~$Q_p$-algebra.
If~$\p_1,\ldots,\p_n$ are the maximal ideals of~$R$ above~$p$ then
\[ R\otimes_Z Z_p \simeq R_{\p_1} \times \ldots \times R_{\p_n}. \]
In particular, each completion~$R_{\p_i}$ is a~$Z_p$-order.

For a maximal ideal~$\p$ of~$R$ and a fractional~$R$-ideal~$I$, we denote by~$I_\p$ the tensor product~$I\otimes_R R_\p$ which is identified with a fractional~$R_\p$-ideal in~$K_\p$.
Given fractional~$R$-ideals~$I$ and~$J$ we have the following equalities:
\[ (I:J)_\p = (I_\p:J_\p),\quad (I+J)_\p = I_\p+J_\p,\quad (IJ)_\p=I_\p J_\p,\quad (I\cap J)_\p = I_\p \cap J_\p.\]

The trace map~$\Tr_{K_\p/Q_p}$ is naturally induced by~$\Tr_{K/Q}$.
Taking completion commutes with taking trace dual, that is,~$(I^t)_\p = (I_\p)^t$.
Hence, the notation~$I^t_\p$ is unambiguous.

A fractional~$R$-ideal~$I$ is called \emph{invertible} if~$I(R:I)=R$.
We say that~$I_\p$ is \emph{locally principal} if~$I_\p$ is isomorphic to~$R_\p$ as an~$R_\p$-module for every maximal ideal~$\p$ of~$R$.
In the next lemma we collect some well known properties.
See for example~\cite[sec.~2]{2022arXiv2MarCMType}.
\begin{lemma}\label{lemma:inv_idl}
    Let~$R$ be an order,~$I$ a fractional~$R$-ideal and~$\p$ a maximal ideal of~$R$.
    Then:
    \begin{enumerate}[(i)]
        \item\label{lemma:inv_idl:loc_princ}~$I$ is invertible if and only if~$I$ is locally principal.
        \item\label{lemma:inv_idl:mult_ring} If~$I$ is invertible then~$(I:I)=R$.
        \item\label{lemma:inv_idl:inverse}~$R \subsetneq (R:\p)$.
        \item\label{lemma:inv_idl:double_inverse}~$ (R:(R:\p))=\p~$.
        \item\label{lemma:inv_idl:p_inv} If~$\p$ is invertible then~$I_\p$ is a principal fractional~$R_\p$-ideal.
        \item\label{lemma:inv_idl:sing} If~$\p$ is not invertible in~$R$ then~$(\p:\p)=(R:\p)$.
        \item\label{lemma:inv_idl:p_cond}~$\p$ is not invertible if and only if~$(R:\cO) \subseteq \p$ if and only if~$R_\p \subsetneq \cO_\p$.
    \end{enumerate}
\end{lemma}

\section{$\p$-equivalence,~$\cS$-equivalence and weak equivalence}
\label{sec:equiv}
In this section, we study several notions of local isomorphism between fractional ideals.

\begin{df}
    Let~$\p$ be a maximal ideal of~$R$.
    We say that two fractional $R$-ideals~$I$ and~$J$ are \emph{$\p$-equivalent} if~$I_\p\simeq J_\p$ as~$R_\p$-modules.
\end{df}

\begin{prop}\label{prop:p-eq}
    Let~$I$ and~$J$ be two fractional~$R$-ideals, and let~$\p$ be a maximal ideal of~$R$.
    Then the following statements are equivalent:
    \begin{enumerate}[(i)]
        \item \label{prop:p-eq:def}~$I$ and~$J$ are~$\p$-equivalent.
        \item \label{prop:p-eq:alpha} There exists a non-zero divisor~$\alpha \in (I:J)_\p$ such that~$I_\p = \alpha J_\p$.
        \item \label{prop:p-eq:loc_colon}~$(I:J)_\p(J:I)_\p=(I:I)_\p$.
        \item \label{prop:p-eq:colon}~$1 \in (I:J)(J:I)+\p$.
    \end{enumerate}
    If any of the above holds, then we have 
    \[ (I:I)_\p=(J:J)_\p,\quad (I:J)_\p = \alpha(I:I)_\p\quad\text{and}\quad(J:I)_\p=(1/\alpha)(I:I)_\p,\]
    with~$\alpha$ as in~\ref{prop:p-eq:alpha}.
\end{prop}
\begin{proof}
    The equivalence of~\ref{prop:p-eq:def} and~\ref{prop:p-eq:alpha} follows from the identifications
    \[ \Hom_{R}(J,I)\otimes_{R} R_\p \simeq (I:J)_\p = (I_\p:J_\p). \]

    It is clear that~\ref{prop:p-eq:alpha} implies~\ref{prop:p-eq:loc_colon}.
    For the converse implication, observe that~$(I:I)_\p$ is a semilocal ring. 
    Hence every invertible fractional~$(I:I)_\p$-ideal is principal.
    By~\ref{prop:p-eq:loc_colon}, we have that~$(I:J)_\p$ is an invertible ideal, hence generated by a non-zero divisor~$\alpha$ as an~$(I:I)_\p$-ideal.
    It follows that
    \[ I_\p = I_\p(I:I)_\p = I_\p \alpha (J:I)_\p \subseteq \alpha J_\p \subseteq I_\p, \]
    which gives~\ref{prop:p-eq:alpha}.

    In the rest of the proof we will use the notation~$N=(I:J)(J:I)$. 
    Now assume~\ref{prop:p-eq:loc_colon}, and
    put~$L = (N+\p) \cap R$.
    We will show that the inclusion~$L \subseteq R$ is everywhere locally an equality.
    Indeed, at~$\p$ we have 
    \[ L_\p = ( N_\p + \p R_\p ) \cap R_\p = (((I:I) + \p) \cap R )_\p = R_\p, \]
    while at~$\q \neq \p$ we have
    \[ L_\q = ( N_\q + \p R_\q ) \cap R_\q = ( N_\q + R_\q ) \cap R_\q = R_\q. \]
    So,~$L=R$, which implies~$1 \in N+\p~$, that is,~\ref{prop:p-eq:colon} holds.
    Conversely, assume~\ref{prop:p-eq:colon}.
    Then~$(I:I) = N +\p(I:I)$.
    Consider the~$R$-module~$M=(I:I)/N$.
    We have
    \[ \frac{M}{\p M} \simeq \frac{(I:I)}{N + \p(I:I)} = 0. \]
    Nakayama's lemma implies that~$M_\p =0$, that is, that~$N_\p = (I:I)_\p$, which is precisely~\ref{prop:p-eq:loc_colon}.

    The final statements are all immediate consequences of~\ref{prop:p-eq:alpha}.
\end{proof}

Let~$R$ be an order and let~$\cS$ be a (possibly infinite) set of maximal ideals of~$R$.
Denote by~$\cS_0$ the subset of~$\cS$ consisting of maximal ideals containing the conductor~$(R:\cO)$ of~$R$.
Note that~$\cS_0$ is a finite set.

\begin{df}
    Two fractional $R$-ideals $I$ and~$J$ are \emph{$\cS$-equivalent} if they are~$\p$-equivalent for every~$\p$ in~$\cS$.
\end{df}

\begin{prop}\label{prop:Seq}
    Let~$I$ and~$J$ be fractional~$R$-ideals. 
    The following statements are equivalent:
    \begin{enumerate}[(i)]
        \item \label{prop:Seq:def}~$I$ and~$J$ are~$\cS$-equivalent.
        \item \label{prop:Seq:colon}~$1 \in (I:J)(J:I) + \p$ for every~$\p$ in~$\cS$.
        \item \label{prop:Seq:S0}~$I$ and~$J$ are~$\cS_0$-equivalent.
        \item \label{prop:Seq:S0colon}~$1 \in (I:J)(J:I) + \p$ for every~$\p$ in~$\cS_0$.
        \item \label{prop:Seq:S0prod}~$1 \in (I:J)(J:I) + \prod_{\p \in \cS_0 }\p$.
    \end{enumerate}
\end{prop}
\begin{proof}
    The equivalences of~\ref{prop:Seq:def} and~\ref{prop:Seq:colon}, and of~\ref{prop:Seq:S0} and~\ref{prop:Seq:S0colon} are consequences of Proposition~\ref{prop:p-eq}.
    If~$\p$ is in~$\cS \setminus \cS_0$ then~$R_\p=\cO_\p$ by Lemma~\ref{lemma:inv_idl}, which implies that~$I$ and~$J$ are~$\p$-equivalent.
    Hence,~\ref{prop:Seq:def} and~\ref{prop:Seq:S0} are equivalent.

    Assume now that~\ref{prop:Seq:S0colon} holds, and set 
    \[ M = \left( (I:J)(J:I) + \prod_{\p \in \cS_0 }\p \right) \cap R. \]
    By assumption,~$M_\p = R_\p$ for every~$\p \in \cS_0$.
    Note that if~$\q$ is a maximal ideal of~$R$ not in~$\cS_0$ then also~$M_\q =R_\q$.
    Therefore~$M=R$ which implies~\ref{prop:Seq:S0prod}.
    The fact that~\ref{prop:Seq:S0prod} implies~\ref{prop:Seq:S0colon} is clear.
\end{proof}

The special case when~$\cS$ consists of all maximal ideals of~$R$ 
gives rise to the notion of weak equivalence.
This was introduced in \cite{DadeTausskyZas}, where a method to effectively compute the weak equivalence classes was also provided.
The algorithm has been improved in \cite{MarsegliaICM} and \cite{2022arXiv2MarCMType}.
\begin{df}
    We say that~$I$ and~$J$ are \emph{weakly equivalent} if they are~$\p$-equivalent for every maximal ideal~$\p$ of~$R$.
\end{df}

\begin{remark}
    The original definition of weak equivalence given in \cite{DadeTausskyZas}, and picked-up in \cite{MarsegliaICM} and \cite{2022arXiv2MarCMType}, uses localization instead of completion.
    The two definitions are equivalent.
\end{remark}

\begin{remark}
    Recall that two fractional~$R$-ideals~$I$ and~$J$ are in the same genus if~$I\otimes_{Z} Z_p \simeq J\otimes_{Z} Z_p$ as~$R\otimes_{Z} Z_p$-modules for every prime of~$Z$.
    As pointed out in \cite[Sec.~5]{LevyWiegand85},~$I$ and~$J$ are weakly equivalent if and only if they are in the same genus.
\end{remark}

\begin{prop}\label{prop:wk-eq}
    Let~$I$ and~$J$ be fractional~$R$-ideals. 
    Denote by~$\frf=(R:\cO)$ the conductor of~$R$.
    The following statements are equivalent:
    \begin{enumerate}[(i)]
        \item \label{prop:wkeq:def}~$I$ and~$J$ are weakly equivalent.
        \item \label{prop:wk:colon}~$1 \in (I:J)(J:I)$.
        \item \label{prop:wk:inv}~$I$ and~$J$ have the same multiplicator ring~$R'=(I:I)=(J:J)$ and there exists an invertible~$R'$-ideal~$L$ such that~$I=LJ$.
        \item \label{prop:wk:peq}~$I$ and~$J$ are~$\p$-equivalent for every maximal ideal~$\p$ containing~$\frf$.
    \end{enumerate}
    If any of the above conditions hold then~$L=(I:J)$.
\end{prop}
\begin{proof}
    The equivalence of~\ref{prop:wkeq:def},~\ref{prop:wk:colon} and~\ref{prop:wk:inv} is \cite[Prop~4.1]{MarsegliaICM} while the last statement is contained in \cite[Cor.~4.5]{MarsegliaICM} for the case~$Z=\Z$.
    The same proofs apply when~$Z$ is any Dedekind domain.
    The equivalence of~\ref{prop:wkeq:def} and~\ref{prop:wk:peq} follows from Proposition~\ref{prop:Seq}.
\end{proof}

We conclude this section by introducing the notation for the sets of classes of the equivalences we have defined so far.
\begin{df}
    Let~$R$ be an order,~$\p$ be a maximal ideal of~$R$ and~$\cS$ be a set of maximal ideals of~$R$. Let:
    \begin{itemize}
        \item~$W(R)$ be the set of weak equivalence classes of fractional~$R$-ideals.
        \item~$\bar W(R)$ be the set of weak equivalence classes of fractional~$R$-ideals~$I$ with~$(I:I)=R$.
        \item~$W_\p(R)$ be the set of~$\p$-equivalence classes of fractional~$R$-ideals.
        \item~$\bar W_\p(R)$ be the set of~$\p$-equivalence classes of fractional~$R$-ideals~$I$ with $(I:I)_\p=R_\p$.
        \item~$W_\cS(R)$ be the set of~$\cS$-equivalence classes of fractional~$R$-ideals.
        \item~$\bar W_\cS(R)$ be the set of~$\cS$-equivalence classes of fractional~$R$-ideals~$I$ with $(I:I)_\p=R_\p$ for every~$\p\in \cS$. 
    \end{itemize}
\end{df}

Since~$\p$-equivalence,~$\cS$-equivalence and weak equivalence are compatible with ideal multiplication, we get that~$W_\p(R)$,~$W_\cS(R)$ and~$W(R)$ inherit the structure of commutative monoids, where the neutral element is the class of~$R$.
Moreover, by Proposition~\ref{prop:wk-eq}, we can partition 
\[ W(R) = \bigsqcup_S \bar W(S), \]
where the disjoint union is taken over all overorders of~$R$.
Similarly, by Proposition~\ref{prop:p-eq}, we have
\[ W_\p(R) = \bigsqcup_{\mathfrak{O}_\p} \bar W_\p(S) \quad\text{and}\quad W_\cS(R) = \bigsqcup_{\mathfrak{O}_\cS} \bar W_\cS(S), \]
where~$\mathfrak{O}_\p$ (resp.~$\mathfrak{O}_\cS$) denotes the set of overorders of~$R$ where we identify two orders if they are locally equal at~$\p$ (resp.~at every~$\p\in\cS$).

\section{Splitting the computation}\
\label{sec:split}
In this section, we describe how the computation of~$\cS$-equivalence classes of a given order~$R$ can be reduced to a computation of weak equivalence classes for a finite set of overorders each having at most one non-invertible maximal ideal.

For every maximal ideal~$\p$ of~$R$, let~$n_\p$ be a positive integer such that~$(\p^{n_\p}\cO)_\p \subseteq R_\p$.
For every~$N\ge n_\p$, we have an equality~$R+\p^{n_\p}\cO = R+\p^N\cO$.
Concretely, one can take~$n_\p = v_p( [\cO:R] )$, where~$p$ is the contraction of~$\p$ in~$Z$ and~$v_p$ is the~$p$-adic valuation. 

\begin{lemma}\label{lemma:RPO}
    Put~$T=R+\p^{n_\p}\cO$. 
    Then:
    \begin{enumerate}[(i)]
        \item \label{lemma:RPO:loc} Locally at~$\p$, we have~$R_\p = T_\p$, while, for every maximal ideal~$\q$ of~$R$ different from~$\p$, we have~$T_\q = \cO_\q$.
        \item\label{lemma:RPO:prime}~$\p T$ is a maximal ideal of~$T$.
        \item\label{lemma:RPO:resfld}~$\p$ and~$\p T$ have isomorphic residue fields~$R/\p \simeq T/\p T$.
        \item\label{lemma:RPO:loc_atp} We have a canonical isomorphism~$T_\p \simeq T_{\p T}$.
        \item\label{lemma:RPO:sing} If~$\p$ is non-invertible then the unique non-invertible maximal ideal of~$T$ is~$\p T$.
    \end{enumerate}
\end{lemma}
\begin{proof}
    Part~\ref{lemma:RPO:loc} is clear from the definition of~$T$.
    Now, consider the~$R$-module $M=T/\p T$. 
    It has finite length, 
    hence it is isomorphic to the direct sum of its localizations.
    Note that~$M_\q = 0$ if~$\q\neq \p$ and~$M_\p = R_\p/\p R_\p \simeq R/\p$.
    This shows that~$\p T$ is a maximal ideal of~$T$ with residue field isomorphic to~$R/\p$, proving~\ref{lemma:RPO:prime} and ~\ref{lemma:RPO:resfld}.

    Note that~$\p T$ is the unique maximal ideal of $T$ containing $\p$.
    To prove~\ref{lemma:RPO:loc_atp}, it is enough to observe that every maximal ideal of~$T$ that meets~$T\setminus \p T$ also meets~$R\setminus \p$.

    Assume now that $\p$ is non-invertible.
    As~$R_\p=T_\p$, it follows that~$(T_\p : \cO_\p) \subsetneq \cO_\p$ by Lemma~\ref{lemma:inv_idl}.\ref{lemma:inv_idl:p_cond}.
    Since~$\p^{n_\p} \subset (T:\cO)$ then~$(T:\cO)$ is a~$\p T$-primary ideal.
    So,~$\p T$ is the unique non-invertible ideal of~$T$, proving~\ref{lemma:RPO:sing}.
\end{proof}

\begin{remark}\label{rmk:type}
    In \cite{2022arXiv2MarCMType}, we showed that the \emph{(Cohen-Macaulay) type} of the order~$R$ at a maximal ideal~$\p$ equals~$\dim_{(R/\p)}(R^t/\p R^t)$.
    If~$T$ is as in Lemma~\ref{lemma:RPO} then~$\typeover{\p}{R} = \typeover{\p T}{T}$. 
    Indeed, the~$R$-module~$N=T^t/\p T^t$ has finite length, and hence it is isomorphic to the direct sum of its localizations.
    Note that~$N_\q = 0$ if~$\q\neq \p$.
    Since~$R_\p =T_\p$, and since completion commutes with taking trace duals, we get~$T^t_\p = R^t_\p$.
    Hence~$ N \simeq N_\p = R^t/\p R^t~$, as required.
\end{remark}

The next two results allow us to compute~$\p$-equivalence and~$\cS$-equivalence classes for an order~$R$ in terms of weak equivalence classes of a finite set of overorders of~$R$, each one having at most one non-invertible maximal ideal. 
\begin{prop}\label{prop:WRWT}
    We have a natural monoid isomorphism
    \[ W_\p(R) \longleftrightarrow W(R+\p^{n_\p}\cO),\]
    which induces a bijection
    \[ \bar W_\p(R) \longleftrightarrow \bar{W}(R+\p^{n_\p}\cO).\]
\end{prop}
\begin{proof}
    Put~$T=R+\p^{n_\p}\cO$. 
    Firstly, observe that every fractional~$T$-ideal is naturally a fractional~$R$-ideal.
    Secondly, two fractional~$T$-ideals~$I$ and~$J$ are weakly equivalent if and only if they are~$\p T$-equivalent, by Lemma~\ref{lemma:RPO}.\ref{lemma:RPO:loc}. 
    By Lemma~\ref{lemma:RPO}.\ref{lemma:RPO:loc_atp} we have a canonical isomorphism~$T_{\p T} \simeq R_\p$,
    which implies that~$I$ and~$J$ are~$\p T$-equivalent if and only if they are~$\p$-equivalent.
    Hence, we have a natural injective map of commutative monoids~$\psi: W(T) \to  W_\p(R)$.
    This map is also surjective since the pre-image of the class in~$W_\p(R)$ of a fractional~$R$-ideal~$I_0$ is given by the class in~$W(T)$ of~$I_0T$.

    Since~$R_\p = T_\p$ by Lemma~\ref{lemma:RPO}.\ref{lemma:RPO:loc}, we deduce that~$\psi(\bar W(T)) \subseteq \bar W_\p(R)$.
    Moreover, if~$I_0$ is a fractional~$R$-ideal with~$(I_0:I_0)_\p = R_\p$, then 
    \[ (I_0T:I_0T)_\p = (I_0 T_\p : I_0 T_\p) = (I_0 R_\p : I_0 R_\p) = (I_0 : I_0)_\p = R_\p=T_\p, \]
    and
    \[ (I_0T:I_0T)_\q = (I_0 T_\q : I_0 T_\q) = (I_0 \cO_\q : I_0 \cO_\q) = \cO_\q = T_\q, \]
    for every maximal ideal~$\q\neq \p$. 
    Hence, the multiplicator ring of~$I_0T$ is~$T$, showing that~$\psi(\bar W(T)) = \bar W_\p(R)$, as required.
\end{proof}

Let~$\cS$ be a set of maximal ideals of~$R$ and denote by~$\cS_0$ the finite subset of~$\cS$ consisting of maximal ideals containing the conductor~$\frf=(R:\cO)$.
The following result, whose proof is similar in spirit to  \cite[Lem.~6.1]{2022arXiv2MarCMType}, allows us to reduce the computation of the monoid of~$\cS$-classes to the computation of a finite number of monoids of weak equivalence classes.

\begin{thm}\label{thm:split}
    We have natural monoid isomorphisms
    \begin{equation}\label{eq:WS_split}
         W_\cS(R) \longleftrightarrow \prod_{\p \in \cS_0} W(R+\p^{n_\p}\cO)
    \end{equation}
    inducing a bijection
    \begin{equation}\label{eq:WSbar_split}
        \bar W_\cS(R) \longleftrightarrow \prod_{\p \in \cS_0} \bar W(R+\p^{n_\p}\cO).
    \end{equation}
\end{thm}
\begin{proof}
    By Proposition~\ref{prop:Seq}, we have that~$W_\cS(R) = W_{\cS_0}(R)$.
    For a maximal ideal~$\p$ in~$\cS_0$ and a fractional~$R$-ideal~$I$, we denote by~$[I]_{\cS_0}$ the class of~$I$ in~$W_{\cS_0}(R)$ and by~$[I]_\p$ the class of~$I$ in~$W_\p(R)$.
    Consider the natural injective map of monoids
    \[ \psi : W_{\cS_0}(R) \longrightarrow \prod_{\p \in \cS_0} W_\p(R). \]
    We now prove that~$\psi$ is surjective.
    Denote by~$\p_1,\ldots,\p_n$ the maximal ideals in~$\cS_0$.
    Consider a vector of classes~$([I_1]_{\p_1},\ldots,[I_n]_{\p_n})$ in the codomain of~$\psi$.
    Without loss of generality, we can assume that the representatives~$I_i$ are contained in~$R$.
    For each~$i$, let~$k_i$ be non-negative integers such that~$\p_i^{k_i}R_{\p_i} \subseteq (I_i)_{\p_i}$, which exist since~$R_{\p_i}$ is Noetherian.
    Put
    \begin{equation}\label{eq:def_J}
       J = \sum_{i=1}^n\left( (I_i+\p_i^{k_i})\prod_{j\neq i} \p_j^{k_j}  \right). 
    \end{equation}
    A direct verification shows that~$J_{\p_i}=(I_i)_{\p_i}$ for each~$i$, and~$J_\q=R_\q$ for every other maximal ideal~$\q$.
    This means that~$[J]_{\cS_0}$ is a preimage of~$([I_1]_{\p_1},\ldots,[I_n]_{\p_n})$, and hence that~$\psi$ is also surjective.
    The isomorphism in Equation~\eqref{eq:WS_split} is the composition of $\Psi$ with the isomorphism from Proposition~\ref{prop:WRWT}.

    For the bijection in Equation~\eqref{eq:WSbar_split}, let~$T=(J:J)$ and assume that for every index~$1\leq i \leq n$ we have~$(I_i:I_i)=R+\p_i^{n_{\p_i}}\cO$.
    Then, by Lemma~\ref{lemma:RPO}, we have
    \[ T_{\p_i} = (I_i:I_i)_{\p_i} = R_{\p_i} \]
    for every~$1\leq i \leq n$. 
    Since we also have that~$T_\q = R_\q$ for every~$\q \not\in\cS_0$ we conclude that~$T=R$.
    Hence, $[J]_{\cS_0}$ is in $\bar W_{\cS_0}(R)$.
\end{proof}

\section{Multiplicator rings and extension}
\label{sec:mult_ring}
Given an inclusion of orders~$R\subseteq T$, Proposition~\ref{prop:presc_mult_ring} below gives conditions on~$T$ that guarantee that the extension to~$T$ of each fractional~$R$-ideal with multiplicator ring~$R$ has multiplicator ring precisely~$T$.
In Proposition~\ref{prop:mult_ring_P_is_T}, we show that the multiplicator ring of any maximal ideal~$\p$ of~$R$ satisfies the equivalent conditions of Proposition~\ref{prop:presc_mult_ring}.

\begin{prop}\label{prop:presc_mult_ring}
    Let~$R$ be an order.
    For an overorder~$R\subseteq T$, the following are equivalent:
    \begin{enumerate}[(i)]
        \item \label{prop:presc_mult_ring:Rt} 
           ~$R^t T$ has multiplicator ring~$T$.
        \item \label{prop:presc_mult_ring:I} 
            For every fractional~$R$-ideal~$I$ with~$(I:I)=R$, the extension~$I T$ has multiplicator ring~$T$.
        \item \label{prop:presc_mult_ring:double_colon_eq}
           ~$(R:(R:T))=T$.
        \item \label{prop:presc_mult_ring:double_colon_incl}
           ~$(R:(R:T))\subseteq T$.
    \end{enumerate}
\end{prop}
\begin{proof}
    Assume~\ref{prop:presc_mult_ring:Rt} and pick~$I$ such that~$(I:I)=R$.
    Note that~$T$ is contained in the multiplicator ring of~$IT$, which in turn is contained in the multiplicator ring of~$IT\cdot I^tT=R^tT$.
    Hence,~$IT$ has multiplicator ring~$T$ as well, which shows~\ref{prop:presc_mult_ring:I}.
    The converse is also true since~$R^t$ has multiplicator ring~$R$.
    Therefore~\ref{prop:presc_mult_ring:Rt} and~\ref{prop:presc_mult_ring:I} are equivalent.

    By Lemma~\ref{lemma:colon_dual} we have
    \[ (R:(R:T))=(R^t(R:T))^t=(R^t(R^tT)^t)^t=(R^tT(R^tT)^t)^t. \]
    So,~\ref{prop:presc_mult_ring:double_colon_eq} holds if and only if~$R^tT(R^tT)^t = T^t$ which is equivalent to~\ref{prop:presc_mult_ring:Rt}, again by Lemma~\ref{lemma:colon_dual}.

    Conditions~\ref{prop:presc_mult_ring:double_colon_eq} and~\ref{prop:presc_mult_ring:double_colon_incl} are equivalent because~$T\subseteq (R:(R:T))$ holds by definition.
\end{proof}

\begin{prop}\label{prop:mult_ring_P_is_T}
    Let~$R$ be an order. 
    For every maximal ideal~$\p$ of~$R$ the multiplicator ring~$T=(\p:\p)$ satisfies the equivalent conditions of Proposition~\ref{prop:presc_mult_ring}.
\end{prop}
\begin{proof}
    If~$\p$ is invertible then~$(\p:\p)=R$ by Lemma~\ref{lemma:inv_idl}.\ref{lemma:inv_idl:mult_ring}, and we are done.
    Assume~$\p$ is not invertible.
    Then, by Lemma~\ref{lemma:inv_idl}.\ref{lemma:inv_idl:double_inverse} and Lemma~\ref{lemma:inv_idl}.\ref{lemma:inv_idl:sing}, we get 
    \[ (R : (R:T)) = (R:(R : (R:\p))) = (R:\p) = T, \]
    as required.
\end{proof}

We conclude the section with Proposition~\ref{prop:mult_ring_P_in_T}, which is not used elsewhere in the paper.
\begin{prop}\label{prop:mult_ring_P_in_T}
    Let~$R$ be an order and~$T$ an overorder of~$R$ satisfying the equivalent conditions of Proposition~\ref{prop:presc_mult_ring}. 
    Then for every maximal ideal~$\p$ of~$R$ containing the conductor~$(R:T)$, we have~$(\p:\p)\subseteq T$.
\end{prop}
\begin{proof}
    If~$T$ is~$R$ then~$(R:T)=R$ which is not contained in any maximal ideal of~$R$.
    So we assume that~$R\subsetneq T$. 
    Note that~$(R:\cO) \subseteq (R:T)$.
    Hence every maximal ideal~$\p$ above~$(R:T)$ is non-invertible by Lemma~\ref{lemma:inv_idl}.\ref{lemma:inv_idl:p_cond}.
    Hence by Lemma~\ref{lemma:inv_idl}.\ref{lemma:inv_idl:sing} and the assumption on~$T$, we have
    \[ (\p:\p)=(R:\p)\subseteq (R:(R:T))=T. \]
\end{proof}

\section{Weak equivalence classes and extension}
\label{sec:wk_ext}
This section explains how, given orders~$R \subset T$, every weak equivalence class of fractional~$R$-ideals with multiplicator ring~$R$ admits a representative contained between two fractional ideals uniquely determined by~$T$.

\begin{lemma}\label{lem:extension}
    Let~$\mathcal{S}$ be a set of maximal ideals of~$R$,
   ~$I_0$ and~$J$ be fractional~$R$-ideals and 
   ~$T$ be an overorder of~$R$.
    Assume that~$I_0T$ is~$\cS$-equivalent to~$J$.
    Then there exists a fractional~$R$-ideal~$I$ which is weakly equivalent to~$I_0$ and satisfies~$IT=J$.
\end{lemma}
\begin{proof}
    Put~$M = (J:I_0T)$.
    After possibly multiplying~$I_0$ by a non-zero divisor of~$K$, we may and do assume that~$M\subseteq R$.
    By Proposition~\ref{prop:p-eq}, for every~$\p \in \cS$, there exists a non-zero divisor~$\alpha_\p \in M_\p$ such that~$\alpha_\p(I_0T)_\p=J_\p$.
    Now, let~$\cS_0$ the intersection of~$\cS$ with the finite set of maximal ideals of~$R$ at which~$I_0T$ and~$J$ are not locally equal, that is, the support of the finite length~$R$-module $(I_0T+J)/(I_0T \cap J)$.
    For each~$\p\in\cS_0$, let~$n_\p$ be a positive integer such that $\p^{n_\p}M_\p \subseteq \alpha_\p R_\p$.
    Define~$\pi$ as the surjective~$R$-linear morphism
    \[ 
        M \twoheadrightarrow  \frac{M}{\left( \prod_{\p \in \cS_0}\p^{n_\p} \right) M } \simeq 
        \prod_{\p \in \cS_0}\frac{M}{\p^{n_\p} M}\simeq 
        \prod_{\p \in \cS_0}\frac{M_\p}{\p^{n_\p} M_\p}.
    \]
    Let~$\alpha$ be any preimage via~$\pi$ of the vector~$(\bar \alpha_\p)_{\p\in \cS_0}$, where~$\bar\alpha_\p$ is the image of~$\alpha_\p$ in~$M_\p/\p^{n_\p} M_\p$.
    Note that~$(\alpha R)_\p = \alpha_\p R_\p$ for every~$\p\in\cS_0$.

    Using the same construction as in the proof of Theorem~\ref{thm:split} (see Equation~\eqref{eq:def_J}), we produce a fractional~$R$-ideal~$L$ satisfying~$L_\p = (\alpha R)_\p$ for every~$\p \in \cS_0$ and~$L_\q=R_\q$ for every~$\q \not \in \cS_0$.
    Set~$I=I_0L$.
    Since~$L$ is a locally principal fractional~$R$-ideal, it is invertible. 
    Hence, by Proposition~\ref{prop:wk-eq},~$I$ and~$I_0$ are weakly equivalent.
    Moreover, for maximal ideals~$\p \in \cS_0$ and~$\q \not \in \cS_0$, we have
    \[ (IT)_\p = (I_0T L)_\p = \alpha_\p I_{0,\p}T_\p = \alpha_\p (I_{0}T)_\p = J_\p, \]
    and
    \[ (IT)_\q = I_{0,\q}L_\q T_\q = I_{0,\q}R_\q T_\q = (I_0T)_\q = J_\q, \]
    which implies that~$IT=J$.
\end{proof}

\begin{prop}\label{prop:good_rep_gen}
    Let~$R\subseteq T$ be orders such that~$R^tT$ has multiplicator ring~$T$.
    Put~$\frf=(R:T)$.
    Let~$J_1,\ldots, J_n$ be representatives of~$\bar W(T)$.
    Let~$I_0$ be a fractional~$R$-ideal with~$(I_0:I_0)=R$.
    Then there exist a unique index~$i$ and a fractional~$R$-ideal~$I$ such that
    \begin{itemize}
        \item~$I$ is weakly equivalent to~$I_0$,
        \item~$IT=J_i$, and
        \item~$\frf J_i \subseteq I \subseteq J_i$.
    \end{itemize}
\end{prop}
\begin{proof}
    By Lemma~\ref{prop:presc_mult_ring} we have that~$I_0T$ has multiplicator ring~$T$.
    By Lemma~\ref{lem:extension}, there exists a fractional~$R$-ideal~$I$ weakly equivalent to~$I_0$ such that~$IT=J_i$, for a unique~$J_i$.
    Now,
    \[ \frf J_i \subseteq \frf IT \subseteq \frf I \subseteq I \subseteq J_i, \]
    as required.
\end{proof}

Proposition~\ref{prop:good_rep_gen} implies that we can find all representatives of the weak equivalence classes with multiplicator ring~$R$ by lifting the sub-$R$-modules of~$J_i/\frf J_i$ for every index~$i=1,\ldots,n$.
We now apply Proposition~\ref{prop:good_rep_gen} to the order~$T=(\p:\p)$.

\begin{thm}\label{thm:good_rep}
    Let~$R$ be an order and~$\p$ be a maximal ideal of~$R$.
    Put~$T=(\p:\p)$.
    Let~$J_1,\ldots,J_n$ be representatives of~$\bar W(T)$.
    Let~$I_0$ be any fractional ideal with multiplicator ring~$R$.
    Then there exist a fractional $R$-ideal~$I$ and a unique index~$i$ such that
    \begin{itemize}
        \item~$I$ and~$I_0$ are weakly equivalent,
        \item~$IT=J_i$, and
        \item~$\p I = \p J_i \subseteq I \subseteq J_i$.
    \end{itemize}
    Moreover, for such~$I$ we have that the natural inclusion~$I\subseteq J_i$ induces an inclusion
    \[ \frac{I}{\p I} \hookrightarrow \frac{J_i}{\p J_i} \]
    of~$R/\p$-vector spaces.
\end{thm}
\begin{proof}
    If~$\p$ is invertible, then~$R=T$ by Lemma~\ref{lemma:inv_idl}.\ref{lemma:inv_idl:mult_ring}, and the statement follows.
    So we assume that~$\p$ is non-invertible.
    Put~$\frf=(R:T)$.
    Note that~$\frf=(R:(R:\p)) = \p$ by Lemma~\ref{lemma:inv_idl}.\ref{lemma:inv_idl:double_inverse}, and that 
   ~$R^tT$ has multiplicator ring~$T$ by Proposition~\ref{prop:mult_ring_P_is_T}.
    Now the statement is an application of Proposition~\ref{prop:good_rep_gen}, together with the observation that~$I\subseteq IT=J_i$ implies~$\p I = \p J_i$.
\end{proof}

\begin{cor}\label{cor:good_rep_inv_ext}
    Let $R$, $\p$, $T$ and $I_0$ be as in Theorem~\ref{thm:good_rep}.
    Assume that~$R^tT$ is invertible as a fractional~$T$-ideal.
    Then there exist a fractional~$R$-ideal~$I$ and a unique index~$i$ such that
    \begin{itemize}
        \item~$I$ and~$I_0$ are weakly equivalent,
        \item~$IT=T$, and
        \item~$\p I = \p \subseteq I \subseteq T$.
    \end{itemize}
    Moreover, for such~$I$ we have that the natural inclusion~$I\subseteq T$ induces an inclusion
    \[ \frac{I}{\p I} \hookrightarrow \frac{T}{\p} \]
    of~$R/\p$-vector spaces.
\end{cor}
\begin{proof}
    Since~$I_0^tI_0=R^t$ by Lemma~\ref{lemma:colon_dual}, we get that~$I_0T$ is also invertible as a fractional~$T$-ideal, that is, $I_0T$ is weakly equivalent to~$T$.
    Hence the result follows from Theorem~\ref{thm:good_rep}, by considering only the invertible class of~$\bar{W}(T)$ and representing it by~$T$.
\end{proof}

We conclude the section with some lemmas that will be used in Section~\ref{sec:wk_lin_alg}.
\begin{lemma}
    \label{lem:wk_inv}
    Let~$R\subseteq T$ be orders and let~$I$ and~$I'$ be fractional~$R$-ideals.
    Assume that~$IT$ has multiplicator ring~$T$.
    Then the following statements are equivalent:
    \begin{enumerate}[(i)]
        \item \label{lem:wk_inv:same_ext}~$I$ and~$I'$ are weakly equivalent, and~$IT=I'T$.
        \item \label{lem:wk_inv:inv_triv} There exists a unique invertible fractional~$R$-ideal~$L$ such that~$I=I'L$ and~$LT=T$.
    \end{enumerate}
\end{lemma}
\begin{proof}
    The fact that~\ref{lem:wk_inv:inv_triv} implies~\ref{lem:wk_inv:same_ext} is clear.
    Now assume~\ref{lem:wk_inv:same_ext}.
    By Proposition~\ref{prop:wk-eq}, there exists a unique invertible fractional~$R$-ideal~$L$ such that~$I=I'L$.
    By extending the equality to~$T$ and using that~$IT=I'T$, we get~$IT=I'LT = IT\cdot LT$.
    This implies that~$LT \subseteq (IT:IT)=T$.
    We also have~$L^{-1}I=I'$.
    By an analogous argument, we obtain~$L^{-1}T \subseteq T$.
    Hence,~$T= L^{-1}T\cdot LT \subseteq T$, that is,~$LT=T$, showing~\ref{lem:wk_inv:inv_triv}.
\end{proof}

\begin{lemma}
    \label{lemma:inv_triv}
    Let~$R\subseteq T$ be orders.
    Consider the function
    \begin{align*} 
        \Phi:\set{\parbox{3.2cm}{\centering invertible fractional\\~$R$-ideals~$L$ such\\ that~$LT=T$}}
        & \longrightarrow
        \prod_\p \frac{T^\times_\p}{R^\times_\p}\\
        L & \longmapsto (x_\p R^\times_\p)_\p
    \end{align*}
    where~$\p$ runs over the maximal ideals of~$R$ and, for each~$\p$,~$x_\p$ is a local generator of~$L$ at~$\p$, that is,~$L_\p = x_\p R_\p$.
    Then~$\Phi$ is a group isomorphism.
\end{lemma}
\begin{proof}
    Let~$L$ be an invertible fractional~$R$-ideal with~$LT=T$.
    Firstly, if~$\p$ is a maximal ideal of~$R$ and~$L_\p = x_\p R_\p$ then~$(LT)_\p=T_\p$ implies that~$x_\p \in T_\p^\times$.
    Secondly, if~$x_\p R_\p = x'_\p R_\p$ then~$x'_\p = u x_\p$ for some unit~$u\in R_\p^\times$, that is, the image~$\Phi(L)$ does not depend on the choice of local generator~$x_\p$.
    Thirdly, observe that~$\Phi$ is a group homomorphism.
    Fourthly, note that~$L$ is in the kernel of~$\Phi$ if and only if, for every maximal ideal~$\p$, the localization~$L_\p$ is generated by a unit in~$R_\p^\times$, that is,~$L_\p = R_\p$.
    Hence~$L=R$, that is,~$\Phi$ is injective.
    Finally, we show that~$\Phi$ is surjective.
    Let~$(x_\p R^\times_\p)_\p$ be an element of the codomain of~$\Phi$.
    For every maximal ideal~$\p$ not containing~$(R:T)$ we have that~$T_\p=R_\p$.
    This implies that, for all but finitely many maximal ideals~$\p$, we have~$x_\p\in R_\p^\times$.
    Using an argument like in the proof of Theorem~\ref{thm:split} (cf.~Equation~\ref{eq:def_J}), we construct a fractional~$R$-ideal~$L$ satisfying~$L_\p = x_\p R_\p$ for every~$\p \supseteq (R:T)$ and~$L_\q = R_\q$ for every other maximal ideal~$\q$.
    Such an ideal~$L$ is locally principal, hence invertible, and since for every~$\p$ we have~$x_\p \in T_\p^\times$, we get~$LT=T$. 
\end{proof}

\section{Weak equivalence classes using linear algebra over residue fields}\
\label{sec:wk_lin_alg}
Let~$R$ be an order and~$\p$ a maximal ideal of~$R$ with residue field~$k=R/\p$ and multiplicator ring~$T=(\p:\p)$.
Let~$J_1,\ldots,J_n$ be representatives of the classes in~$\bar W(T)$.
Consider the~$k$-algebra~$A=T/\p$ and the~$k$-vector spaces~$\cJ_i=J_i/\p J_i$.
We have a natural action by multiplication of~$A$ on each~$\cJ_i$ satisfying~$\cJ_iA = \cJ_i$.
For each index~$1\leq i\leq n$, consider the natural bijection
\begin{align*}
        \Psi_i : \set{\parbox{4cm}{\centering fractional~$R$-ideals~$I$\\ such that~$\p J_i \subseteq I \subseteq J_i$}} & \longrightarrow 
\set{\parbox{2.3cm}{\centering sub-$k$-vector\\ spaces of~$\cJ_i$}}\\
    I & \longmapsto I/\p J_i.
\end{align*}

For each index~$i$, let~$\cV_i$ denote the set of sub-$k$-vector spaces~$\cI$ of~$\cJ_i$ satisfying
\begin{itemize}
    \item~$\set{ a \in A : a\cI \subseteq \cI } = k$, and 
    \item~$\cI A=\cJ_i$.
\end{itemize}

\begin{lemma}
    \label{lemma:psi_i}
    The set~$\cV_i$ is the image via~$\Psi_i$ of the set of fractional~$R$-ideals~$I$ satisfying~$\p J_i \subseteq I \subseteq J_i$,~$(I:I)=R$ and~$IT=J_i$.
\end{lemma}
\begin{proof}
    Let~$I$ be a fractional~$R$-ideal in the domain of~$\Psi_i$ and put~$\cI=\Psi_i(I)$.
    Note that~$IT=J_i$ if and only if~$\cI A=\cJ_i$.
    If so, then
    \begin{equation}\label{eq:mult_ring_p}
        \set{a \in A : a\cI\subseteq \cI} = (T \cap (I:I))/\p = (I:I)/\p,
    \end{equation}
    where the second equality follows from~$T=(J_i:J_i)=(IT:IT)\supseteq (I:I)$.
    Hence, Equation~\eqref{eq:mult_ring_p} gives that~$\set{a \in A : a\cI\subseteq \cI} = k$ if and only if~$(I:I)=R$.
\end{proof}

\begin{thm}\label{thm:wk_classes_only_linear_algebra}\
    Consider the group~$U=\prod_\p T_\p^\times/R_\p^\times$, where~$\p$ runs over the set of maximal ideals of~$R$. The group~$U$ acts freely on each~$\cV_i$, and the collection of bijections~$\Psi_i$ induces a bijection
    \[ \Psi:\bar W(R) \longrightarrow \bigsqcup_{i=1}^n\frac{\cV_i}{U}. \]
\end{thm}
\begin{proof}
    By Lemma~\ref{lemma:inv_triv}, the group~$U$ is isomorphic to the group of invertible~$R$-ideals~$L$ satisfying~$LT=T$, which acts by multiplication on the set of fractional~$R$-ideals~$I$ such that~$(I:I)=R$ and~$IT=J_i$ for each index~$i$.
    Hence, we obtain the desired action of~$U$ on~$\cV_i$ by Lemma~\ref{lemma:psi_i}.
    Note that the action is free by Lemma~\ref{lem:wk_inv}.

    By Theorem~\ref{thm:good_rep}, for every fractional~$R$-ideal~$I_0$ with~$(I_0:I_0)=R$ there is an index~$i$ and a weakly equivalent ideal~$I$ such that~$IT=J_i$. 
    This condition is equivalent to have~$\Psi_i(I) \in \cV_i$ by Lemma~\ref{lemma:psi_i}.
    It follows that~$\Psi$ is surjective.
    
    Now, let~$I'$ be another fractional~$R$-ideal with multiplicator ring~$R$.
    A necessary condition for~$I$ and~$I'$ to be weakly equivalent is that~$IT$ and~$I'T$ are both weakly equivalent to~$J_i$, which we assume.
    Hence, by Lemma~\ref{lem:extension}, we can and do assume that~$IT=I'T=J_i$, that is,~$I$ and~$I'$ are both in the domain of the same map~$\Psi_i$. 
    So~$\Psi_i(I)$ and~$\Psi_i(I')$ are both in~$\cV_i$.
    Lemmas~\ref{lem:wk_inv} and~\ref{lemma:inv_triv} tell that~$I$ and~$I'$ are weakly equivalent if and only if~$\Psi_i(I)$ and~$\Psi_i(I')$ are in the same~$U$-orbit.
    Therefore~$\Psi$ is also injective.
\end{proof}

The results contained in Section~\ref{sec:split} allow us to reduce the computation of weak equivalence classes of a given order to the computation of the weak equivalence classes of a finite number of overorders each having a unique non-invertible maximal ideal.
The next proposition can be used to describe the group~$U$ from Theorem~\ref{thm:wk_classes_only_linear_algebra} when~$R$ has a unique non-invertible maximal ideal~$\p$.

\begin{prop}\label{prop:U_one_sing}\
    Let~$R \subsetneq S$ be orders.
    Let~$\frf$ be a fractional~$S$-ideal such that~$\frf\subsetneq R$.
    Assume that~$R$ has a unique non-invertible maximal ideal~$\p$ and that the only maximal ideal of~$R$ containing~$\frf$ is~$\p$.
    Then we have group isomorphisms
    \[ \frac{(S/\frf)^\times}{(R/\frf)^\times} 
    \overset{\vphi_1}{\simeq}
    \prod_\q \frac{S_\q^\times}{R_\q^\times}
    \overset{\vphi_2}{\simeq}
    \frac{S_\p^\times}{R_\p^\times}, \]
    where the product is taken over all maximal ideals~$\q$ of~$R$.
\end{prop}
\begin{proof}
    We have that~$(R:\cO) \subseteq (R:S) \subsetneq R$, where the second inclusion is strict because~$R\subsetneq S$.
    By Lemma~\ref{lemma:inv_idl},~$\p$ is the unique maximal ideal of~$R$ containing the conductor~$(R:\cO)$.
    Hence,~$\p$ is the unique maximal ideal containing~$(R:S)$.
    The existence of the isomorphism $\vphi_2$ in the statement is a consequence of the fact that~$R_\q=S_\q$ for every maximal ideal~$\q$ of~$R$ not containing~$(R:S)$.
    The isomorphism $\vphi_1$ is constructed in the same way as in \cite[Prop.~12.11]{Neukirch} where the author treats the special case when $K$ is a number field, $S=\cO$ and $\frf$ is the conductor of $R$.
\end{proof}
\begin{remark}\label{rmk:U_effective}
    Assume that the order~$R$ has a unique non-invertible maximal ideal~$\p$.
    Applying Proposition~\ref{prop:U_one_sing} to the inclusion of orders~$R\subsetneq T=(\p:\p)$ gives us at least two methods to compute the group~$U$ and its action on each set of subspaces~$\cV_i$ as in Theorem~\ref{thm:wk_classes_only_linear_algebra}.

    The first one is obtained by setting~$\frf=(R:T)$, which equals~$\p$ by Lemma~\ref{lemma:inv_idl}, in Proposition~\ref{prop:U_one_sing}.
    In this case, we get 
    \[ U \simeq \frac{A^\times}{k^\times}. \]
    If~$k$ is a finite field, then also~$A$ is a finite ring, giving an efficient method to compute~$U$.
    This is the case when~$Z=\Z$ or~$Z=\F_q[t]$, for example.
    Since the set of representatives of~$U$ so obtained is in~$A$, it is easy to compute the orbits of the action of~$U$ on each~$\cV_i$. 

    The second method to compute~$U$ is obtained by setting~$\frf=(R:\cO)$, the conductor of~$R$, in Proposition~\ref{prop:U_one_sing}.
    Then we get
    \[
        U \simeq \frac{(T/\frf)^\times}{(R/\frf)^\times},
    \]
    which is a quotient of a subgroup of the ray class group~$(\cO/\frf)^\times$.
    Hence, if~$Z=\Z$ or~$Z=\F_q[t]$ then we can compute~$U$ by applying the algorithms described in \cite{klupau05}.
    
    Both methods above allow us to represent~$U$ by elements in~$T$.
    Hence, we can realize the action of~$U$ on each set~$\cV_i$ by acting by multiplication on the set of ideals~$\Psi_i^{-1}(\cV_i)$.
\end{remark}

\section{Algorithms}
\label{sec:alg}
We now turn the results presented in the previous sections into two algorithms, \ref{alg:WRbar} and~\ref{alg:WR}, to compute~$\bar{W}(R)$ and~$W(R)$ for an order~$R$, respectively.
Both algorithms require enumerating sub-vector spaces of a finite dimensional vector spaces over the residue field of a (non-invertible) maximal ideal~$\p$ of~$R$, see Remarks~\ref{rmk:alg:WRbar}.\eqref{rmk:alg:WRbar:enum} and~\ref{rmk:alg:both}.\eqref{rmk:alg:both:enum}.
Hence to guarantee that the enumeration is finite we need to assume that~$R/\p$ is finite.
Under this assumption, we obtain that the quotient of every two $Z$-lattices is finite.
In particular, the index $[\cO:R]$ is finite.

The implementation in Magma \cite{Magma}, for the case~$Z=\Z$, is part of the package \textrm{AlgEtQ} available at \url{https://github.com/stmar89/AlgEt}.

\begin{alg}[\customlabel{alg:WRbar}{\textup{\textbf{ComputeWbar}}}]\
    \newline\noindent
    Input: an order~$R$.\\
    Output: a set~$\bar{\cW}$ of representative of the classes of~$\bar{W}(R)$.\\
    \begin{enumerate}[label=({\arabic*}), ref={\arabic*},wide=0pt, leftmargin=\parindent]\vspace{-1em}
        \item \label{alg:WRbar:early_exit} if~$R = \cO$ then \textbf{return}~$\bar\cW = \set{R}$;
        \item \label{alg:WRbar:sing} compute the set~$\cS_0=\set{\p_1,\ldots,\p_n}$ of non-invertible maximal ideals of~$R$;
        \item \label{alg:WRbar:ngt1} if~$n>1$ then
        \begin{enumerate}[label=({\ref{alg:WRbar:ngt1}.\alph*}), ref={\ref{alg:WRbar:ngt1}.\alph*}]
            \item \label{alg:WRbar:ngt1:for} for each~$1\leq i \leq n$ do
            \begin{enumerate}[label=({\ref{alg:WRbar:ngt1:for}.\roman*}), ref={\ref{alg:WRbar:ngt1:for}.\roman*}]
                \item compute the contraction~$p_i$ of~$\p_i$ in~$Z$;
                \item set~$n_i=v_{p_i}([\cO:R])$;
                \item recursively use~\ref{alg:WRbar} to compute a set of representatives~$\bar{\cW}_i$ of~$\bar{W}(R+\p_i^{n_i}\cO)$;
            \end{enumerate}
            \item use Theorem~\ref{thm:split} (and the construction described in the proof) to compute~$\bar\cW$ from~$\prod_{i=1}^n \bar\cW_i$;
            \item \textbf{return}~$\bar\cW$;
        \end{enumerate}
        \item \label{alg:WRbar:neq1} if~$n = 1$ then
        \begin{enumerate}[label=({\ref{alg:WRbar:neq1}.\alph*}), ref={\ref{alg:WRbar:neq1}.\alph*}]
            \item set~$\p=\p_1$;
            \item compute~$T=(\p:\p)$,~$k=R/\p$ and~$A=T/\p$; 
            \item \label{alg:WRbar:U} compute a set of representatives~$U=T_\p^\times/R_\p^\times$ using Remark~\ref{rmk:U_effective}
            \item if~$R^tT$ is invertible in~$T$, then set~$m=1$ and~$J_1=T$; else
            recursively use~\ref{alg:WRbar} to compute a set of representatives~$\set{J_1,\ldots,J_m}$ of ~$\bar{W}(T)$;
            \item \label{alg:WRbar:neq1:for} for each~$1\leq j\leq m$ do
            \begin{enumerate}[label=({\ref{alg:WRbar:neq1:for}.\roman*}), ref={\ref{alg:WRbar:neq1:for}.\roman*}]
                \item set~$\cJ_j=J_j/\p J_j$;
                \item \label{alg:WRbar:neq1:Vj} compute the set~$\cV_j$ of sub-$k$-vector spaces~$\cI$ of~$\cJ_j$ satisfying~$\cI A=\cJ_i$ and~$\set{a \in A : a\cI \subseteq \cI}=k$;
            \end{enumerate}
            \item \label{alg:WRbar:neq1:Wbar} compute~$\bar{\cW}$ from~$\sqcup_{j=1}^m \cV_j/U$ (see Remark~\ref{rmk:alg:WRbar}.\eqref{rmk:alg:WRbar:orbit});
            \item \textbf{return}~$\bar{\cW}$.
        \end{enumerate}
    \end{enumerate}
\end{alg}
\begin{thm}\label{thm:alg:WRbar}
    Algorithm~\ref{alg:WRbar} is correct.
\end{thm}
\begin{proof}
    If~$R$ is maximal, that is, there are no non-invertible maximal ideals, then every fractional~$R$-ideal is invertible, that is, weakly equivalent to~$R$ by Proposition~\ref{prop:wk-eq}.

    If~$R$ is non-maximal then the recursion is finite because the orders~$R+\p_i^{n_{\p_i}}\cO$ and~$(\p_i:\p_i)$ are both strictly bigger than~$R$ by Lemmas~\ref{lemma:RPO} and~\ref{lemma:inv_idl}, respectively.
    The output is a set of representatives of~$\bar{W}(R)$ by Theorem~\ref{thm:split} when~$n>1$, and by Theorem~\ref{thm:wk_classes_only_linear_algebra} when~$n=1$, together with Corollary~\ref{cor:good_rep_inv_ext} in the case that~$R^tT$ is an invertible fractional~$T$-ideal.
\end{proof}
\begin{remark}\label{rmk:alg:WRbar}
    This remark contains a list of comments about various steps of Algorithm~\ref{alg:WRbar}.
    \begin{enumerate}[wide=0pt, leftmargin=1.5\parindent]
        \item \label{rmk:alg:WRbar:type} The early-exit call in Step \eqref{alg:WRbar:early_exit} can be generalized.
        Recall the definition of local (Cohen-Macaulay) type of an order~$R$
        from Remark~\ref{rmk:type}.
        If the maximum local type of~$R$ over the non-invertible maximal ideals is~$1$ then~$R$ is also said to be \emph{Gorenstein} and~$\bar{W}$ consists only of the class of~$R$; see for example \cite[Prop.~3.4]{2022arXiv2MarCMType}.
        We have a complete description of~$\bar{W}(R)$ also when the maximum of the local types of~$R$ is~$2$.
        In this case, by \cite[Thm.~6.2]{2022arXiv2MarCMType}, a fractional~$R$-ideal with multiplicator ring~$R$ is, locally at each maximal ideal~$\p$, isomorphic to~$R$ or~$R^t$. 
        Hence using \cite[Lem.~6.1]{2022arXiv2MarCMType}, we can effectively reconstruct~$\bar{W}(R)$.
        \item \label{rmk:alg:WRbar:sing} In Step \eqref{alg:WRbar:sing}, we need to compute the non-invertible maximal ideals of~$R$.
        This can be achieved by first computing the conductor~$\frf=(R:\cO)$ of~$R$ and then lift the maximal ideals of the Artinian ring~$R/\frf$.
        Observe that under our assumption that~$R$ is residually finite, we have that~$R/\frf$ is a finite ring.
        \item \label{rmk:alg:WRbar:enum} In Step \eqref{alg:WRbar:neq1:Vj}, we need to find all sub-$k$-spaces~$\cI$ of the~$A$-module~$\cJ_i$ that satisfy two conditions:~$\set{a\in A : a\cI\subseteq \cI}=k$ and~$\cI A=\cJ_i$.

        We start by discussing the second. 
        Let~$\cK_1,\ldots,\cK_r$ be the maximal sub-$A$-modules of~$\cJ_j$. 
        Observe that~$\cI A=\cJ_j$ if and only if~$\cI \not\subset \cK_i$ for all~$1\leq i\leq r$.
        Also, if~$\cI\subseteq \cK_i$, for some~$i$, then all sub-$k$-spaces of~$\cI$ will also be contained in~$\cK_i$.
        Hence, it is wise to look for the subspaces~$\cI$ of~$\cJ_j$ by recursively enumerating maximal sub-$k$-spaces, and stop the recursion if we find a sub-$k$-space contained in one of the~$\cK_i$.
        Observe that if~$\cJ_i = T/\p$ then each~$\cK_i$ is of the form~$\mathfrak{P}_i/\p$ where~$\mathfrak{P}_i$ is a maximal ideal of~$T$ above~$\p$.

        We now consider the first condition,~$\set{a\in A : a\cI\subseteq \cI}=k$.
        Let~$\frB'$ be a minimal set of generators of~$A$ over~$k$.
        Assume that~$1\in\frB'$ and set~$\frB=\set{a \in \frB' : a\neq 1}$.
        Let~$\cI$ be a sub-$k$-space of~$\cJ_j$, with~$k$ basis~$\mathfrak{I}$.
        Then,~$\set{a\in A : a\cI\subseteq \cI}=k$ if and only if for all~$a\in \frB$ there exists~$x\in \mathfrak{I}$ such that~$ax \not\in \cI$.
        One can also proceed as follows. 
        For each~$a\in A$, denote by~$\mathcal{M}_a$ the set of minimal sub-$k[a]$-modules of~$\cJ_j$. 
        Set~$\mathcal{M} = \cup_{a\in\mathfrak{B}}\mathcal{M}_a$.
        Observe that~$a\cI\subseteq \cI$ implies that there exists~$M\in \mathcal{M}_a\subseteq \mathcal{M}$ such that~$M\subseteq \cI$.
        If for every~$M\in\mathcal{M}$ we have~$M\not\subset \cI$ then we have the equality~$\set{a\in A : a\cI\subseteq \cI}=k$.
        From this perspective, and in contradiction with the previous one, it is preferable to look for the subspaces~$\cI$ of~$\cJ_j$ by recursively enumerating minimal sub-$k$-spaces, and stop the recursion if we find a sub-$k$-vector space containing some~$M \in \mathcal{M}$.
        \item \label{rmk:alg:WRbar:orbit} In Step~\eqref{alg:WRbar:neq1:Wbar}, we need to extract a set of representative~$\bar{\mathcal{W}}$ for the disjoint union of the sets~$\cV_j$.
        For each~$j$, initialize an empty set~$\mathfrak{O}$. 
        Then, whenever the enumeration of the sub-$k$-spaces of~$\cJ_j$ produces an~$\cI$ in~$\cV_j$ such that~$\cI$ is not in~$\mathfrak{O}$, then add the coset~$\cI U$ to~$\mathfrak{O}$ and a lift of~$\cI$ via the quotient map~$J_j\to \cJ_j$ to~$\bar{\cW}$. 

        Another way to do it, which allows us to skip the computation of~$U$ in Step \eqref{alg:WRbar:U}, is the following: for each~$j$, whenever the enumeration of the sub-$k$-spaces of~$\cJ_j$ produces an~$\cI$ in~$\cV_j$, compute its lift~$I$ via the quotient map~$J_j\to \cJ_j$; if there is no ideal in~$\bar{\cW}$ which is weakly equivalent to~$I$, which can be checked by using Proposition~\ref{prop:wk-eq}, then append~$I$ to~$\bar{\cW}$.

        The difference between the two approaches is that, in the first, we need to check equality of vector spaces, while, in the second, we have to check weak equivalence of fractional ideals, which can be done in polynomial time; see Section~\ref{sec:running_time_comparison}.
    \end{enumerate}
\end{remark}

\begin{alg}[\customlabel{alg:WR}{\textup{\textbf{ComputeW}}}]\
    \newline\noindent
    Input: an order~$R$.\\
    Output: a set~$\cW$ of representative of the classes of~$W(R)$.\\
    \begin{enumerate}[wide=0pt, leftmargin=\parindent]\vspace{-1em}
        \item if~$R = \cO$ then \textbf{return}~$\cW = \set{R}$;
        \item compute the set~$\cS_0=\set{\p_1,\ldots,\p_n}$ of non-invertible maximal ideals of~$R$;
        \item \label{alg:WR:ngt1} if~$n>1$ then
        \begin{enumerate}[label=({\ref{alg:WR:ngt1}.\alph*}), ref={\ref{alg:WR:ngt1}.\alph*}]
            \item \label{alg:WR:ngt1:for} for each~$1\leq i \leq n$ do
            \begin{enumerate}[label=({\ref{alg:WR:ngt1:for}.\roman*}), ref={\ref{alg:WR:ngt1:for}.\roman*}]
                \item compute the contraction~$p_i$ of~$\p_i$ in~$Z$;
                \item set~$n_i=v_{p_i}([\cO:R])$;
                \item recursively use~\ref{alg:WR} to compute a set of representatives~$\cW_i$ of~$W(R+\p_i^{n_i}\cO)$;
            \end{enumerate}
            \item use the construction described in the proof of Theorem~\ref{thm:split} to compute~$\cW$ from~$\prod_{i=1}^n \cW_i$;
            \item \textbf{return}~$\cW$;
        \end{enumerate}
        \item \label{alg:WR:neq1} if~$n = 1$ then
        \begin{enumerate}[label=({\ref{alg:WR:neq1}.\alph*}), ref={\ref{alg:WR:neq1}.\alph*}]
            \item \label{alg:WR:neq1:Wbar} use~\ref{alg:WRbar} to compute a set of representatives~$\bar{\cW}_R$ of~$\bar W(R)$;
            \item \label{alg:WR:n1:minoo} Compute the minimal overorders~$S_1,\ldots,S_m$ of~$R$;
            \item \label{alg:WR:neq1:rec} for each~$1\leq j \leq m$, recursively use~\ref{alg:WR} to compute a set of representatives~$\cW_{S_j}$ of~$W(S_j)$ while keeping track of the multiplicator rings (that is, use the decomposition~$W(S_j) = \bigsqcup \bar{W}(T)$, where~$T$ runs over the overorders of~$S_j$, to partition~$\cW_{S_j} = \bigsqcup \bar{\cW}_T$);
            \item \label{alg:WR:neq1:disjunion} set~$\cW = \bar{\cW}_R \sqcup \bigsqcup_j \bar{\cW}_{S_j} \sqcup \bigsqcup_T \bar{\cW}_T$, where~$T$ runs over the strict overorders of all~$S_1,\ldots,S_m$, avoiding repetitions.
            \item \textbf{return}~$\cW$.
        \end{enumerate}
    \end{enumerate}
\end{alg}

\begin{thm}\label{thm:alg:WR}
    Algorithm~\ref{alg:WR} is correct.
\end{thm}
\begin{proof}
    If~$R$ is maximal, that is, there are no non-invertible maximal ideal, then every fractional~$R$-ideal is invertible, that is, weakly equivalent to~$R$ by Proposition~\ref{prop:wk-eq}.
    The over orders~$R+\p_i^{n_i}$ are strictly bigger than~$R$ by Lemma~\ref{lemma:RPO}. 
    The same is true for the minimal overorders~$S_1,\ldots,S_m$ of R.
    Therefore the recursions terminate.
    The output contains all representatives of classes in~$\bar{W}(R)$ by Theorem~\ref{thm:split} when~$n>1$ and by construction for~$n=1$.
\end{proof}

\begin{remark}\label{rmk:alg:both}\
    This remark contains a list of comments about various steps of Algorithm~\ref{alg:WR}.
    \begin{enumerate}[wide=0pt, leftmargin=1.5\parindent]
        \item Observe that Remark~\ref{rmk:alg:WRbar}.\eqref{rmk:alg:WRbar:sing} is relevant also for Algorithm~\ref{alg:WR}.
        \item \label{rmk:alg:both:enum} In Step \eqref{alg:WR:n1:minoo} in Algorithm~\ref{alg:WR}, we need to compute the minimal overorders of an order~$R$ having a unique non-invertible maximal ideal~$\p$.
        One can affectively compute them using \cite[Rmk.~5.4]{HofmannSircana20}.
        \item In Step \eqref{alg:WR:n1:minoo} in Algorithm~\ref{alg:WR}, the minimal overorders have common overorders. 
        For example,~$S_j\subseteq (\p:\p)$ for each~$j$, by \cite[Prop.~5.2]{HofmannSircana20}.
        So, in the recursions in both Algorithms, in order to avoid recomputing the same information multiple times, it is important to cache the lattice of overorders of each order~$T$ together with~$\bar{W}(T)$.
        \item Algorithms~\ref{alg:WRbar} and~\ref{alg:WR} can easily be modified to compute~$\bar{W}_\cS(R)$ and~$W_\cS(R)$ for some set of maximal ideals~$\cS$ of~$R$ by using Theorem~\ref{thm:split}.
    \end{enumerate}
\end{remark}

\section{Running time and comparison with previous algorithm}
\label{sec:running_time_comparison}

\subsection*{Running time of~\ref{alg:WRbar} and~\ref{alg:WR}}
In this section, we assume that~$K$ is an \'etale algebra over~$\Q$, that is, a finite product of number fields, and that~$R$ is a~$\Z$-order in~$K$.
Once a~$\Q$-basis of~$K$ is fixed, orders and fractional ideals can be described using~$\Z$-basis and multiplication matrices.
In what follows, we will use this as input, and the running times will be expressed with respect to the size of this input.
Given fractional~$R$-ideals~$I$ and~$J$, one can compute~$I+J$,~$I\cap J$,~$IJ$,~$(I:J)$, and~$J/I$ (if~$I\subseteq J$) in polynomial time, see for example \cite[Prop.~5.8]{BuchmannLenstra94}.

In Step \eqref{alg:WRbar:sing} in both algorithms, we need to compute the non-invertible maximal ideals of the order~$R$.
This can be done by lifting to~$R$ the maximal ideals of the finite ring~$R/dR$ where~$d$ is the discriminant of the order~$R$, and for each such lift~$\p$ check whether it is not invertible; cf.~\cite[Rmk~5.8]{HofmannSircana20}.
Alternatively, one could pre-compute the maximal order~$\cO_K$, factor~$d\cO_K$, take intersections with~$R$, and check for non-invertibility.
The second approach has the advantage that the information over~$\cO_K$ can be cached and re-used in the recursion of the algorithms.
In both cases, the resulting procedure is subexponential.

Note that the reduction from the case of an order with more than one non-invertible maximal ideal to the case of an order with a single non-invertible maximal ideal, that is, the application of Theorem~\ref{thm:split} in Step \eqref{alg:WRbar:ngt1} in both algorithms, is done in polynomial time.

We now start to analyze Step \eqref{alg:WRbar:neq1} of Algorithm~\ref{alg:WRbar}, from which we borrow the notation.
The computation is trivial in the cases discussed in Remark~\ref{rmk:alg:WRbar}.\eqref{rmk:alg:WRbar:type}.
We will assume that we are not in such cases.


The most expensive part of the algorithm is Step \eqref{alg:WRbar:neq1:Vj}.
In the enumeration of the sub-vector spaces of~$\cJ_j$, we can take advantage of early exits in the recursive process, as explained in Remark~\ref{rmk:alg:WRbar}.\eqref{rmk:alg:WRbar:enum}.
These conditions seem hard to control, so we will pretend that~$\cV_j$ consists of all sub-vector spaces of~$\cJ_j$.
The number of such spaces is a function of the size~$q$ of the finite field~$R/\p$ and of the dimension of~$\cJ_j$, which is bounded from above by the dimension~$d$ of~$\cO/\p \cO$; see \cite[Cor.~2.2]{Greither82}.
In particular, we get that the number of such subspaces grows exponentially with~$d$.
As pointed out in Remark~\ref{rmk:alg:WRbar}.\eqref{rmk:alg:WRbar:orbit}, for each sub-vector space we have to check if it is already on the keep-track list~$\mathfrak{O}$, and if not to compute its orbit under the action of~$U$.
The latter task depends on the size of~$U$, while the former requires a very large number of equality checks.

Now we discuss Step \eqref{alg:WR:neq1} of Algorithm~\ref{alg:WR}.
The computation of the minimal overorders in Step \eqref{alg:WR:n1:minoo} shares many similarities with Step \eqref{alg:WR:neq1:Wbar}. 
Indeed, by \cite[Prop~5.2]{HofmannSircana20}, a minimal overorder~$S$ of the order~$R$ with unique non-invertible maximal ideal~$\p$ satisfies~$S\subseteq (\p:\p)$. 
This implies that the minimal overorders of~$R$ can be found by enumerating the sub-$R/\p$-spaces of~$(\p:\p)/R$.
Finally, Step \eqref{alg:WR:neq1:disjunion} is trivial if one keeps track of the multiplicator rings in Step~\eqref{alg:WR:neq1:rec}.
\subsection*{Comparison with previous algorithms}
In the seminal work \cite{DadeTausskyZas}, the authors introduce the notion of weak equivalence classes for an order~$R$ and give a method to compute~$\cW(R)$.
They show that every class in~$W(R)$ admits a representative~$J$ satisfying~$\frf \subseteq J \subseteq \cO$, where~$\frf=(R:\cO)$ is the conductor of~$R$; see \cite[Sec.~2.3]{DadeTausskyZas}.
Hence, one can compute~$W(R)$ by enumerating sub-$R$-modules of the finite quotient~$\cO/\frf$, and sieve from their preimages via the quotient map~$\cO\to \cO/\frf$ a minimal set of representatives.

An improvement on this method has been described by the author in \cite{MarsegliaICM}.
This approach uses the decomposition
\[ W(R) = \bigsqcup_S \bar W(S), \]
where~$S$ runs over the overorders of~$R$.
Note that the overorders can be computed using \cite{HofmannSircana20}.
We are left to compute each~$\bar W(S)$.
In \cite[Prop.~5.1]{MarsegliaICM}, it is shown that every class in~$\bar W(S)$ admits a representative~$J$ such that~$(S:T) \subseteq J \subseteq T$, where~$T$ is any overorder of~$S$ such that the extension~$S^tT$ of the trace dual ideal~$S^t$ of~$S$ is invertible as a fractional~$T$-ideal.
This means that it is sufficient to enumerate sub-$S$-modules of the quotient~$T/(S:T)$ which is potentially (much) smaller than~$\cO/\frf$ since we have inclusions
\[ \frf \subseteq (S:T) \subseteq S \subseteq T \subseteq \cO. \]
This is a big gain in efficiency, even if we have to perform an enumeration process for every overorder.

The gist of Theorem~\ref{thm:wk_classes_only_linear_algebra} and of Algorithm~\ref{alg:WRbar} is that we can reduce the enumeration process by considering the even smaller quotients~$J/\p J$ where~$J$ runs over representatives of~$\bar\cW((\p:\p))$.
On the one hand there are orders where the new approach does not provide any advantage; see Example~\ref{ex:bad}.
On the other hand, Example~\ref{ex:good} and the description of the computation of abelian varieties given in Section~\ref{sec:example_app} shows that practically the gain in performance is striking.

\section{Examples and applications}
\label{sec:example_app}
\begin{example}\label{ex:good}
    Consider the algebra~$K=\Q[x]/h$ where
    \[ h=(x^2+6x+25)(x^4-6x^3+42x^2-150x+625). \]
    Write~$K=K_1\times K_2$, with~$K_1$ a quadratic number field and~$K_2$ a quartic one.
    Denote by~$\pi$ the class of~$x$ in~$K$ and consider the involution defined by~$\bar{\pi} = 25/\pi$.
    Let~$\cO_K$,~$\cO_{K_1}$ and~$\cO_{K_2}$ be the corresponding maximal orders.
    The maximal ideals above~$2$ of~$\cO_K$ are
    \[ \frP_0=2\cO_{K_1} \times \cO_{K_2}, \quad \frP_1=\cO_{K_1} \times \p_1, \quad \frP_2=\cO_{K_1} \times \p_2, \quad \frP_3=\cO_{K_1} \times \p_3 \]
    with~$\p_1 = \bar{\p}_1$ and~$\p_2 = \bar{\p}_3$ maximal ideals of~$\cO_{K_2}$.

    Put 
    \[\frf = \frP_0^7\frP_1^2\frP_2^2\frP_3^4 \quad\text{and}\quad R=\Z[\pi] + \frf.\]
    One can check that~$R$ is an order with conductor~$(R:\cO_K)=\frf$.
    If one would like to compute~$\bar W(R)$ using \cite{DadeTausskyZas} one would have to enumerate the~$303348$ sub-$R$-modules of~$\cO_K/\frf$.

    Among the overorders~$S$ of~$R$ such that~$R^t S$ is invertible as a fractional~$S$-ideals, there is a unique one which is minimal with respect to inclusion.
    We denote it by~$S_0$.
    Computing~$\bar W(R)$ using \cite{MarsegliaICM} would require enumerating through the~$98990$
    sub-$R$-modules of~$S_0/(R:S_0)$.

    One can compute that~$R$ has a unique non-invertible maximal ideal ideal~$\p$.
    Let~$T=(\p:\p)$.
    We compute that~$\bar W(T)$ consists of~$46$ weak equivalence classes, say represented by~$J_1,\ldots, J_{46}$.
    Algorithm~\ref{alg:WRbar} requires us to enumerate the sub-$R/\p$-vector spaces of~$J_i/\p J_i$, for each~$i$:~$16$ of them have~$67$ subspaces, while the remaining~$30$ have~$16$ subspaces.
\end{example}

\begin{example}\label{ex:bad}
    Let~$K$ be an \'etale algebra with maximal order~$\cO_K$.
    Let~$R$ an order in~$K$ whose conductor~$\p=(R:\cO_K)$ is a maximal ideal.
    Then the~\ref{alg:WRbar} as well as the algorithms built out of \cite{DadeTausskyZas} and \cite{MarsegliaICM} described in Section~\ref{sec:running_time_comparison} require looping through the~$R$-submodules of~$\cO_K/\p$.
    In other words, for such an order $R$, the enumeration steps in all three algorithms coincide.
    Note that it is easy to construct order whose conductor is a maximal ideal.
    For example, take~$R=Z+p\cO_K$, where~$p$ is any prime of~$Z$.
\end{example}

On the one hand, Example~\ref{ex:bad} shows that there is no considerable performance gain in using the algorithm presented in this paper for certain orders.
On the other hand, a practical application described in the rest of this section shows that the performance gain is striking for most orders.

\subsection*{Ideal Classes}
Let~$R$ be an order.
We denote by~$\ICM(R)$ the set of isomorphism classes of fractional~$R$-ideals.
Ideal multiplication endows~$\ICM(R)$ of the structure of a commutative monoid.
Analogous to the partition~$ W(R) = \bigsqcup_S \bar W(S)$, with~$S$ running over the overorders of~$R$, we also have 
\[ \ICM(R) = \bigsqcup_S \bar{\ICM}(S), \]
where~$\bar\ICM(S)$ denotes the set of isomorphism classes with multiplicator ring~$S$.
For an overorder~$S$, let~$\Pic(S)$ denote the \emph{Picard group} of~$S$, that is, the group of invertible fractional~$S$-ideals up to isomorphism.
As shown in \cite[Thm.~4.6]{MarsegliaICM} and \cite[section 5]{2022arXiv2MarCMType}, we have that~$\Pic(S)$ acts freely on~$\bar\ICM(S)$ and the quotient of this action is precisely~$\bar W(S)$.
Therefore, we obtain an algorithm to compute~$\ICM(R)$ by combining the results contained in this paper, together with the results contained in \cite{klupau05} which give an effective method to compute~$\Pic(S)$. 

\subsection*{Abelian varieties over finite fields}
Let~$h(x)$ be a Weil polynomial determining, by Honda-Tate theory, an isogeny class~$\cI_h$ of abelian varieties over the finite field~$\F_q$ of dimension~$g$.
Assume that~$h(x)$ is squarefree, that is, has no repeated complex roots, and that~$q$ is a prime number or~$h(x)$ is ordinary, that is, the coefficient of~$x^g$ is coprime with~$q$.
Set~$K=\Q[x]/h$ and~$\pi = x\bmod h(x)$.
By \cite[Cor~4.4]{MarsegliaAbVar}, we have a bijection
\[ \set{\parbox{3.8cm}{\centering abelian varieties in~$\cI_h$\\up to~$\F_q$-isomorphism}} \longleftrightarrow \ICM(\Z[\pi,q/\pi]).  \]

By using this bijection, we computed the isomorphism classes of abelian varieties over finite fields for~$615269$ isogeny classes over various dimensions~$g$ and finite fields~$\F_q$.
The computation of the weak equivalence classes for all orders of all the isogeny classes was achieved in~$70$ CPU-hours
on an Intel Xeon CPU E5-2687W v2 3.40GHz using no more than $1$ GB of RAM (only $11$ isogeny classes required more than $300$ MB of RAM).
The computation was performed using the implementation of~\ref{alg:WR} and~\ref{alg:WRbar} which is available at \url{https://github.com/stmar89/AlgEt}\footnote{See the branch $\mathrm{recursive\_wk\_icm}$, commit $1831318$.}.
A previous run of the same computation on the same machine using an implementation of the results contained in \cite{MarsegliaICM} was not able to finish in more than one hundred CPU-months.
This shows that the current algorithm provides a performance boost of several orders of magnitude.

A subset of the output of this computation, together with several other interesting arithmetic information, has been uploaded in a joint effort with Edgar Costa, Taylor Dupuy, David Roe, Christelle Vincent and Mackenzie West to the LMFDB~\cite{lmfdb} (currently available through \url{https://abvar.lmfdb.xyz/Variety/Abelian/Fq/}).

\bibliographystyle{amsalpha}
\renewcommand{\bibname}{References} 
\providecommand{\bysame}{\leavevmode\hbox to3em{\hrulefill}\thinspace}
\providecommand{\MR}{\relax\ifhmode\unskip\space\fi MR }
\providecommand{\MRhref}[2]{%
  \href{http://www.ams.org/mathscinet-getitem?mr=#1}{#2}
}
\providecommand{\href}[2]{#2}

\end{document}